\documentclass[preprint, 5p, times, sort&compress]{elsarticle}

\usepackage{microtype}

\usepackage{graphicx}
\usepackage{amsmath}
\usepackage{amssymb}
\usepackage{bm}
\usepackage{mathtools}

\usepackage{booktabs}
\usepackage{siunitx}
\usepackage{array}
\usepackage{colortbl}
\usepackage{algorithm}
\usepackage{algpseudocode}
\usepackage{tcolorbox}

\usepackage[htt]{hyphenat}
\usepackage{hyperref}

\def\vecr{\bm{r}}
\def\veck{\bm{k}}
\def\veca{\bm{a}}
\def\vecb{\bm{b}}
\def\vecc{\bm{c}}
\def\vecy{\bm{y}}
\def\vtau{\bm{\tau}}
\def\vrho{\bm{\rho}}
\def\veta{\bm{\eta}}
\def\zbot{z_\text{e}}
\def\dd{\mathrm{d}}
\def\imag{\mathrm{i}}
\def\ee{\mathrm{e}}
\def\sd{\mathord{:}}
\def\vzero{\bm{0}}

\begin{document}
        \begin{frontmatter}
            \title{A fast and accurate computation method for reflective diffraction simulations\tnoteref{tag:title}}
            \author[uec]{Shuhei Kudo\corref{cor1}}
            \ead{shuhei-kudo@uec.ac.jp}
            \author[uec]{Yusaku Yamamoto}
            \ead{yusaku.yamamoto@uec.ac.jp}
            \author[tottori]{Takeo Hoshi}
            \ead{hoshi@tottori-u.ac.jp}
            \affiliation[uec]{organization={The University of Electro-Communications},
                addressline={1-5-1 Chofugaoka},
                city={Chofu},
                postcode={182-8585},
                country={Japan}
            }
            \affiliation[tottori]{organization={Tottori University},
                addressline={4-101 Koyama-cho Minami},
                city={Tottori},
                postcode={680-8550},
                country={Japan}}
            \cortext[cor1]{Corresponding author}
            
            \begin{abstract}
                We present a new computation method for simulating
                reflection high-energy electron diffraction
                and the total-reflection high-energy positron diffraction experiments. The two experiments are used commonly for the structural analysis of material surface. 
                The present paper improves the conventional numerical method, the multi-slice method, for faster computation, since the present method avoids the matrix-eigenvalue solver for the computation of matrix exponentials and can adopt higher-order ordinary differential equation solvers. 
                Moreover, we propose a high-performance implementation based on multi-thread parallelization and cache-reusable subroutines.
                In our tests, this new method performs up to 2,000 times faster than the conventional method.
            \end{abstract}            
        \end{frontmatter}

        \section{Introduction}
        Nowadays, fast and accurate numerical methods are of great importance for data-driven science, in particular, for a global search analysis. As a typical problem, an equation is characterized by a parameter set $\bm{x}=(x^{(1)},x^{(2)},..,x^{(n)})$ and one should solve the equation numerically on many different points in the parameter space of $\bm{x}$. A major application field in computational physics is the reverse analysis of experimental measurement,
        in which the observed data $\bm{d}$ can be described as a function of the target quantity $\bm{x}$
        ($\bm{d}_\mathrm{cal}=f(\bm{x})$).
        Hereinafter, the function $f$ is called the \emph{forward} model.
        A typical approach is an optimization analysis
        in which the optimal value of the target quantity
        $\bm{x}^* = \operatorname*{argmin}_{\bm{x}} | \bm{d} - F(\bm{x}) |^2$ is computed by a global optimization method,
        such as a grid search or Bayesian optimization.
        Another typical approach is Bayesian inference, in which
        the posterior priority density $P(\bm{x}|\bm{d})$ is obtained as a histogram.
        These approaches require the computation of $f(\bm{x})$
        for a large dataset of $\bm{x}=\bm{x}_1, \bm{x}_2, \ldots$.
        Among such approaches,
        a rapid numerical computation is desired for performing
        reverse analysis with a larger dataset.

        The present paper is motivated by the reverse analysis of the two experimental measurements,
        reflection high-energy electron diffraction (RHEED)~\cite{Ichimiya2004} and
        total-reflection high-energy positron diffraction (TRHEPD)~\cite{Fukaya2019a,HUGENSCHMIDT_2016_SurfSciRep_rev}.
        RHEED and TRHEPD are experimental probes for crystal surface structures,
        i.e.,\ positions of the atoms $\bm{x}$ in the surface and subsurface atomic layers.
        In RHEED and TRHEPD, quantum beams of  electrons and positrons, respectively, are irradiated to the crystal surface and observe
        the diffraction patterns of the reflected waves $\bm{d}$.
        TRHEPD is more sensitive to the shallower layers of the surface than RHEED.
        There is an open-source RHEED\slash TRHEPD simulation software \texttt{sim-trhepd-rheed}~\cite{Hanada2022, str-github},
        which is used for 
        many data analysis of surface structures based of RHEED~\cite{RHEED_SURF_SCI_1993, RHEED_HIKITA_JVST_1993, RHEED-Kudo-SurfIntf-1994, RHEED-HANADA-SurfSci-1994, HANADA_PRB_1995,RHEED-YAMADA-PRL-1995,RHEED-OHTAKE-PRB-1999, RHEED-Miotto-ApplPhysLett-1999, RHEED-Ohtake-JCrysGrow-1999, RHEED-Ohtake-PRB60-1999, RHEED-RealTime-Ohtake-PRB60-1999, RHEED-Miotto-2000,RHEED-Ohtake-PRB-2001,RHEED-Ohtake-PRB-2002}
        and TRHEPD~\cite{TANAKA2020_ACTA_PHYS_POLO,  HOSHI2022_TRHEPD_SION,TANAKA_2023_JJAP,Motoyama2022,TSUJIKAWA_2022_MOLECULES, TSUJIKAWA_2022_PRB}.
        Motoyama et.al.~\cite{Motoyama2022} demonstrated the reverse analysis of TRHEPD
        using their fully automated software 2DMAT
        that enable us the global search algorithms.
        In several global search algorithms, the forward problems with many points of the parameter space of $\bm{x}=\bm{x}_1, \bm{x}_2, \ldots$ are solved simultaneously as a parallel computation. 

        The present standard method for performing a forward computation of RHEED\slash TRHEPD is the \emph{multi-slice method}~\cite{Ichimiya1983,Ichimiya2004}.
        The multi-slice method solves the stationary Schr\"{o}dinger equation as
        a boundary value problem (BVP) by expanding the wave function with periodic functions
        in the x-y plane and obtaining ordinary differential equations (ODEs) for the z-axis.
        The range of $z$ is divided into thin \emph{slice}s of size $h$
        and the potential $v(z)$ is approximated by a constant matrix $A_i$ in each slice.
        Thus, the matrix exponential $\hat{A}_i=\ee^{hA_i}$ becomes the transfer matrix of each slice $i$,
        and their product $\hat{A}_0 \hat{A}_1\cdots$ becomes that of the whole crystal.
        To solve this, Ichimiya invented a technique that we call the \emph{recursive reflection technique}. The technique solves the matrix eigenvalue problem for the calculation of matrix exponential $A_i$ on each slice, which will be costly, when the matrix dimension of $A_i$ increases. The computational cost is proportional to the number of the slices $N_{\rm slice}$ and, thus, to the inverse of slice size $h$. 

        We propose a new fast computation method by reorganizing the problem as a matrix ODE and
        improving the recursive reflection technique in higher-order ODE solver algorithms. 
        The present method realizes fast and accurate numerical computation, because the present method does not require a matrix-eigenvalue problem and the slice size $h$ can be chosen to be  
        more than $10$ fold larger than the original one to achieve the same accuracy.
        Furthermore, we propose a high-performance implementation method based on the implementation techniques in the highe-performance computing (HPC),
        such as multi-thread parallelization and the active use of cache-reusable subroutines,
        to exploit the performance of recent CPUs.
        As a result, our method performs up to $2,000$ times faster than the implementation of the conventional method in our tests.
        Although we apply our new method only to RHEED\slash TRHEPD simulation in this study, this method can potentially be extended to
        other many-beam reflective diffraction simulations.

        \begin{table}
                \centering \caption{List of notations used in this paper}\label{tab:notations}
                \begin{tabular}{ll}
                        \toprule
                        Symbol & Description \\
                        \midrule
                        $\imag$ & the imaginary unit\\
                        $\ee$ & the Euler's number\\
                        $\bm{a}$ & vectors\\
                        $(\bm{a})_i$ & the $i$-th component of the vector $\bm{a}$.\\
                        $A$ & matrices\\
                        $(A)_{i,j}$ & the $(i, j)$ component of the matrix $A$ \\
                        $I_n$ & the identity matrix of size $n\times n$\\
                        $O_n$ & the zero matrix of size $n\times n$ \\
                        $\hat{X}$ & numerical approximations of quantity $X$ \\
                        \bottomrule
                \end{tabular}
        \end{table}

        We summarize the notations used in this article in Table~\ref{tab:notations}.

        The reminder of the article is organized as follows.
        In Section 2, we introduce the simulation model of RHEED\slash TRHEPD
        and summarize the conventional computation technique
        developed by Ichimiya~\cite{Ichimiya1983}.
        In Section 3, we describe our proposed method,
        and we compare the performance of our proposed method with that of the conventional implementation in Section 5.
        We present our discussions and conclusions in Section 6 and 7.


        \section{Simulation technique for RHEED\slash TRHEPD}
        \subsection{Numerical problem}
        \begin{figure}
                \begin{minipage}{0.48\linewidth}
                        \centering\includegraphics[width=0.8\linewidth]{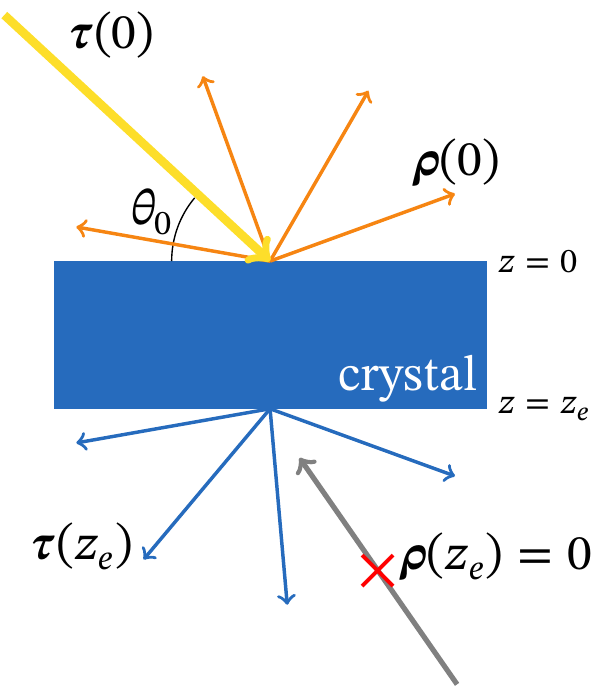}
                \end{minipage}%
                \hspace{0.04\linewidth}%
                \begin{minipage}{0.48\linewidth}
                        \centering\includegraphics[width=0.8\linewidth]{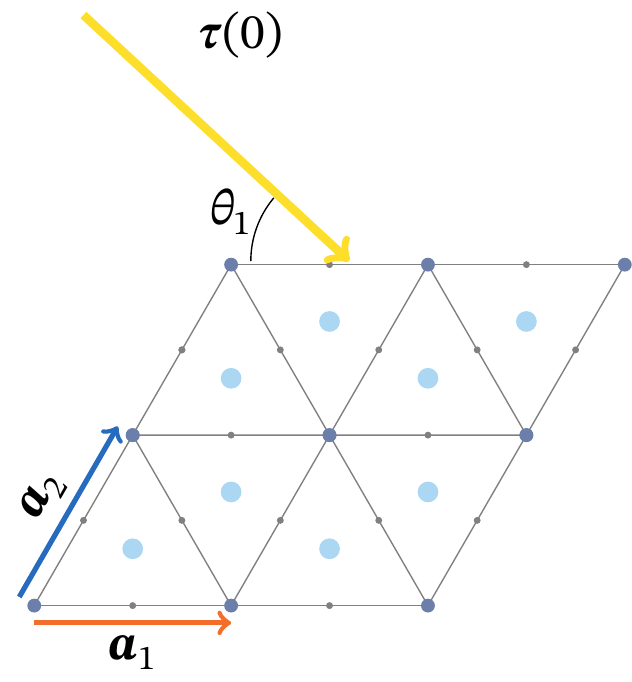}
                \end{minipage}
                \caption{Left: Illustration of the boundary conditions of the computational system, a shallow surface region of crystal.
                        $\theta_0$ is the altitude of the incident wave $\vtau(0)$.
                        Right: Illustration of the lattice structure of the crystal surface ($z=0$) from the top.
                        $\theta_1$ is the azimuth of $\vtau(0)$.} \label{fig:lattice}
        \end{figure}
        Here, we introduce the numerical method for RHEED developped by Ichimiya~\cite{Ichimiya1983,Ichimiya2004} and TRHEPD  ~\cite{Fukaya2019a}.
        The computational system is a shallow surface region of crystal, as shown in Fig.~\ref{fig:lattice}. The coordinate axes are set so that the material surface is parallel to the x-y plane at precisely $z=0$, and the position of the bottom boundary of the system is $z=\zbot$.
        The positions of the atoms are two-dimensionally periodic parallel to the x-y plane,
        but not periodic along the z-axis.

        Let us assume that the incident wave and reflection wave satisfy the stationary Schr\"odinger equation:
        \begin{equation}
                \left(\Delta + \gamma^2 + v(\vecr)\right) \psi(\vecr) = 0, \label{eq:sch1}
        \end{equation}
        where $\vecr$ is the position,
        $v$ is the potential generated by the atoms in the system, and
        $\gamma$ is the wave number of electrons or positrons in a vacuum.
        The periodicity of $v$ can be written by the  lattice vectors $\veca_1$ and $\veca_2$ that are parallel to the surface, as shown in the right panel of fig.~\ref{fig:lattice} :
        \begin{equation}
                v(\vecr) = v(\vecr + m_1\veca_1+m_2\veca_2), \quad m_1, m_2 \in \mathbb{Z}. \label{eq:pot1}
        \end{equation}
        Following Bloch's theorem \cite{KITTEL-TEXT,Ichimiya2004}, 
        $\psi$ is a product of a periodic function $\phi$ and the plane wave $\ee^{\imag\vecb\cdot\vecr}$:
        \begin{align}
                \psi(\vecr) &= \ee^{\imag\vecb\cdot\vecr}\phi(\vecr),\\
                \phi(\vecr) &= \phi(\vecr + m_1\veca_1 + m_2\veca_2). \label{eq:bloch2}
        \end{align}

        Assuming that the wave is continuous, we rewrite the equations in the frequency domain,
        which naturally satisfy the periodic boundary conditions:
        \begin{align}
                \hat{v}(\vecr) &= \sum_{j=1}^\infty u_j(z)\ee^{\imag\veck_j\cdot \vecr_{\mathrm{xy}}}, \\
                \hat{\psi}(\vecr) &= \ee^{\imag\vecb_0\cdot \vecr_{\mathrm{xy}}}\sum_{j=1}^\infty c_j(z)\ee^{\imag\veck_j\cdot \vecr_{\mathrm{xy}}},
        \end{align}
        where $\hat{v}$ and $\hat{\psi}$ are the frequency-domain counterparts of $v$ and $\psi$, respectively,
        $u_j$ and $c_j$ are the frequency components of $v$ and $\psi$, respectively,
        $\veck_j$ are the reciprocal rods,
        $\vecr_{\mathrm{xy}}$ is the component of $\vecr$ which is parallel to the x-y plane,
        $z$ is the z component of $\vecr$, and $\vecb_0$ is the projected wavevector of $\vecb$ to the x-y plane.
        To perform the numerical computation, we truncate the series by $n$ components and substitute them into Eq.~\eqref{eq:sch1}: 
        \begin{equation}
                \left(\frac{\dd^2}{\dd z^2}-\left\|\vecb_0+\veck_j\right\|^2+\gamma^2\right)c_j(z) + \sum_{k=1}^n u_{l(j,k)}(z)c_k(z) = 0 \label{eq:sch2}
        \end{equation}
        where $1\le j \le n$ and we introduce an index conversion function $l(j, k)$ calculated from the position of the reciprocal rods.
        The number of components $n$ is $n=10^1 -- 10^2$ in the present paper. 
        We further simplify the equation by defining a vector $\vecc$, matrix $U$, and diagonal matrix $\Gamma$:
        \begin{equation}
                \frac{\dd^2}{\dd z^2}\vecc(z) = -(U(z) + \Gamma^2)\vecc(z), \label{eq:sch3}
        \end{equation}
        where $(\vecc)_j=c_j$, $(U)_{j,k}=u_{l(j,k)}$, and $(\Gamma)_{j,j}=\sqrt{\gamma^2-\|\vecb_0+\veck_j\|^2}$.

        The boundary condition of the problem is Robin (third type).
        Next, we define the incident wave $\vtau$ and reflected wave $\vrho$ as follows:
        \begin{align}
                \vtau(z) &= \left(-\mathrm{i}\frac{\dd}{\dd z}+\Gamma\right)\vecc(z), \label{eq:vtau} \\
                \vrho(z) &= \left(\mathrm{i}\frac{\dd}{\dd z}+\Gamma\right)\vecc(z). \label{eq:vrho}
        \end{align}
        Then, we let the $\veck_1$ component represent a plane wave, which means that $\veck_1=\vzero$, and
        the boundary conditions are defined as follows:
        \begin{align}
                \vtau(0) &= \begin{pmatrix}1&0&\cdots&0\end{pmatrix}^\top, \label{eq:bd0}\\
                \vrho(\zbot) &= \mathbf{0}. \label{eq:bd1}
        \end{align}
        The reflection wave at the top $\vrho(0)$ is the term that we aim to compute, and $\vtau(\zbot)$ is not used.

        The RHEED\slash TRHEPD intensity $\hat{\veta}$ is computed from $\vrho(0)$ as follows:
        \begin{equation}
                (\veta)_i = \begin{cases}
                        \|(\vrho(0))_i\|^2 \frac{\sin\theta_0\sqrt{\gamma^2-\|\vecb_0\|^2}}{\sqrt{\gamma^2-\|\vecb_0+\veck_i\|^2}} &
                        \text{if } \gamma^2 > \|\vecb_0+\veck_i\|^2 \\
                        0 & \text{otherwise}
                \end{cases},
        \end{equation}
        where $\theta_0$ is the altitude of the incident wave.

        Because the measurements are performed for many pairs of the altitudes $\theta_0$ and azimuths $\theta_1$,
        simulations must be performed for each angle as well.
        The intensity $\veta$ as a function of $\theta_0$ and $\theta_1$ is called
        the \emph{rocking curve},
        which is $\bm{d}_\mathrm{cal}$ for RHEED\slash TRHEPD in the reverse analysis.

        \subsection{Conventional method}
        The \emph{multi-slice method} is used as a standard method for the simulation of RHEED\slash TRHEPD.
        The method splits the z-range to make thin \emph{slices} of the system
        and approximates the coefficient matrix function by a step-wise constant function.
        As a result, the whole problem can be written as the product of the transfer matrices for each slice.
        Thus, the problem is mathematically converted into simultaneous linear equations,
        but a special technique is required to solve these equations to avoid numerical breakdowns.

        \paragraph{Discretization scheme}
        Let us rewrite Eq.~\eqref{eq:sch3} with $\vtau$ and $\vrho$
        by using the equations $\vtau+\vrho=2\Gamma\vecc$ and $\vrho-\vtau=2\mathrm{i}\frac{\dd}{\dd z}\vecc$:
        \begin{align}
                \frac{\dd}{\dd z}\vtau(z) &= \frac{\imag}{2}U(z)\Gamma^{-1}(\vtau(z)+\vrho(z)) + \mathrm{i}\Gamma\vtau(z), \\
                \frac{\dd}{\dd z}\vrho(z) &= -\frac{\imag}{2}U(z)\Gamma^{-1}(\vtau(z)+\vrho(z)) - \mathrm{i}\Gamma\vrho(z).
        \end{align}
        We then obtain a block-wise matrix equation:
        \begin{align}
                \frac{\dd}{\dd z}
                \begin{pmatrix}
                        \vtau(z)\\
                        \vrho(z)
                \end{pmatrix}
                &=
                A(z)
                \begin{pmatrix}
                        \vtau(z)\\
                        \vrho(z)
                \end{pmatrix},
                \label{eq:tr1}\\
                A(z) &\coloneqq
                \begin{pmatrix}
                        \frac{\imag}{2}U(z)\Gamma^{-1}+\imag\Gamma & \frac{\imag}{2}U(z)\Gamma^{-1} \\
                        -\frac{\imag}{2}U(z)\Gamma^{-1} & -\frac{\imag}{2}U(z)\Gamma^{-1} -\imag\Gamma
                \end{pmatrix}.
        \end{align}
        The conventional method approximates the coefficient matrix $A(z)$ with a step-wise constant matrix.
        Now, let the crystal be split into $L$-slices at the point $z_0=\zbot, z_1, \ldots, z_L=0$, where $z_i<z_{i+1}$ for $i=0$ to $L-1$.
        We approximate $A(z)$ as follows:
        \begin{align}
                A(z) &\approx A_i \coloneqq A\left(\frac{z_i+z_{i+1}}{2}\right), \quad \text{where} \ \ z_i \le z < z_{i+1}. \label{eq:ds0}
        \end{align}
        In each slice, Eq.~\eqref{eq:tr1} is approximated by a linear ODE with constant coefficient matrices,
        thus, the exact solution for each approximated equation can be written as follows:
        \begin{equation}
                \begin{pmatrix}
                        \hat{\vtau}(z_{i+1}) \\
                        \hat{\vrho}(z_{i+1})
                \end{pmatrix}
                =
                \mathrm{e}^{(z_{i+1}-z_i)A_i}
                \begin{pmatrix}
                        \hat{\vtau}(z_{i}) \\
                        \hat{\vrho}(z_{i})
                \end{pmatrix}.
                \label{eq:ichimiya2}
        \end{equation}
        Hence, by letting $\hat{A}_i \coloneqq \ee^{(z_{i+1}-z_i)A_i}$, we obtain a linear simultaneous equation:
        \begin{equation}
                \begin{pmatrix}
                        \hat{\vtau}(0) \\
                        \hat{\vrho}(0)
                \end{pmatrix}
                =
                \left(\hat{A}_{L-1}\hat{A}_{L-2}\cdots\hat{A}_0\right)
                \begin{pmatrix}
                        \hat{\vtau}(\zbot) \\
                        \hat{\vrho}(\zbot)
                \end{pmatrix}.
                \label{eq:ichimiya3}
        \end{equation}
        This equation is solved by a block approach. Let the product of matrices be split into $n\times n$ sub-matrices:
        $\begin{pmatrix} \hat{X} & \hat{Y} \\
                \hat{Z} & \hat{W}\end{pmatrix}$,
        and we have
        \begin{align}
                \hat{\vtau}(0) &= \hat{X} \hat{\vtau}(\zbot) + \hat{Y}\hat{\vrho}(\zbot), \\
                \hat{\vrho}(0) &= \hat{Z} \hat{\vtau}(\zbot) + \hat{W}\hat{\vrho}(\zbot).
        \end{align}
        By inserting the boundary condition Eq.~\eqref{eq:bd1}, and equating out $\hat{\vtau}(\zbot)$, we find
        \begin{align}
                \hat{\vrho}(0) &= \hat{Z} \hat{X}^{-1} \hat{\vtau}(0).\label{eq:ichimiya4}
        \end{align}

        \paragraph{Recursive reflection technique}
        It is difficult to solve Eq.~\eqref{eq:ichimiya4} because the matrix product has a very large condition number.
        Instead of computing the whole product, the \emph{recursive reflection technique} inductively
        constructs the reflection matrix.
        First, the reflection matrix of the bottom slice is computed, and
        then, this matrix is combined with the matrix for the slice above the bottom to obtain the next reflection matrix.
        By iteratively repeating this process, we obtain the reflection matrix of the whole crystal.

        Let $\hat{A}_i=\begin{pmatrix}
                \hat{X}_i & \hat{Y}_i \\
                \hat{Z}_i & \hat{W}_i
        \end{pmatrix}$ where all sub-matrices are $n\times n$.
        By inserting the boundary condition Eq.~\eqref{eq:bd1} for Eq.~\eqref{eq:ichimiya2} for $i=0$, we have
        \begin{equation}
                \hat{\vrho}(z_1) = \hat{Z}_0 \hat{X}_0^{-1} \hat{\vtau}(z_1).
        \end{equation}
        For $\hat{R}_0 \coloneqq \hat{Z}_0 \hat{X}_0^{-1}$, we have $\hat{\vrho}(z_1) = \hat{R}_0\hat{\vtau}(z_1)$.
        Let us suppose that $\hat{\vrho}(z_{i}) = \hat{R}_{i-1}\hat{\vtau}(z_i)$ for $i>0$; then,
        we have
        \begin{align}
                \hat{\vtau}(z_{i+1}) &= \hat{X}_i\hat{\vtau}(z_i) + \hat{Y}_i \hat{\vrho}(z_i) = \left(\hat{X}_i + \hat{Y}_i\hat{R}_{i-1}\right)\hat{\vtau}(z_i), \\
                \hat{\vrho}(z_{i+1}) &= \hat{Z}_i\hat{\vtau}(z_i) + \hat{W}_i \hat{\vrho}(z_i) = \left(\hat{Z}_i + \hat{W_i}\hat{R}_{i-1}\right)\hat{\vtau}(z_i).
        \end{align}
        Equating out $\hat{\vtau}(z_i)$ from the equation, we obtain the next reflection matrix $\hat{R}_i$ as
        \begin{align}
                \hat{R}_i &\coloneqq  \left(\hat{Z}_i+\hat{W_i}\hat{R}_{i-1}\right)\left(\hat{X}_i+\hat{Y}_i\hat{R}_{i-1}\right)^{-1}, \label{eq:ri1}\\
                \hat{\vrho}(z_{i+1}) &= \hat{R}_i\hat{\vtau}(z_{i+1}).
                \label{eq:ri2}
        \end{align}
        By induction, the desired reflection matrix $\hat{R}_{L-1}=\hat{Z}\hat{X}^{-1}$ can be computed from the bottom to the top.
        Note that we are assuming that the matrix inverses always exist.
        Then, the boundary condition problem becomes $\hat{\vrho}(0) = \hat{R}_{L-1}\hat{\vtau}(0)$, a matrix-vector multiplication.

        \paragraph{Treatment for the bulk layer}
        In RHEED/TRHEPD experiments, the target structure (atom positions) is located at the very shallow surface region, whereas
        the reminder of the system is called the \emph{bulk layer} in which the atom positions are fixed to be that in the ideal (known) crystals.
        The conventional method assumes that the reflection matrix $R_i$ of the bulk layer rapidly converges to $R_\mathrm{bulk}$.
        Thus, the conventional method stops the computation if $R_i$ converges
        and skips the rest of the computation for the bulk layer. It is noted $R_\mathrm{bulk}$ is computed before the reverse analysis and reused in the structure search.
        We are using same technique in our proposed method.

        \paragraph{Strategy for a faster computational method}
        Here one can find that the conventional method is based on an unbalanced strategy, because it perform a high-cost numerical procedure (eigenvalue problem solver) that rises from  
        low-order approximations. Therefore, a faster numerical method can rise, if one reorganizes the strategy. 
        The origin of the high-cost procedure in the conventional method is based on the computation of matrix exponentials.  
        There are several methods for computing matrix exponentials, including eigenvalue decomposition-based methods
        and the ``scaling-and-squaring'' method, but all such methods require $2n\times 2n$ matrix computations~\cite{Moler2003}.
        Ichimiya proposed a technique to reduce the matrix size
        from $2n\times 2n$ to $n\times n$ by using the structure of the matrix,
        but this approach is still costly.
        Moreover, the step-wise constant coefficient scheme used in the conventional method is
        second-order; thus, a small slice size $h$ or a large number of the slices $L$ is required to perform an accurate computation.
        The discretization scheme used in the conventional method is known as
        the second-order Magnus method~\cite{Blanes2008},
        which has interesting features such as structure preservation.
        The present numerical problem, however, does not have such a structure. 

        \section{Proposed method}
        We propose a new technique for solving the BVP Eq.~\eqref{eq:sch3} based on two new steps.
        In the first step, we rewrite the equation as a matrix ODE, which has a simpler form and thus can be more easily applied to well-known ODE solvers, such as the Runge-Kutta method.
        In the second step, we rewrite the recursive reflection technique as a linear transformation by post-multiplication of a matrix.
        As a result, we can apply the recursive reflection technique to a wide variety of ODE solvers.
        We also propose a fast implementation technique based on the knowledge of HPC.

        \subsection{Matrix ODE}
        We rewrite Eq.~\eqref{eq:sch3} by using $\vecc$ and $\frac{\dd}{\dd z}\vecc$ instead of $\vtau$ and $\vrho$.
        Let $\vecy(z) = \begin{pmatrix}
                        \vecc(z) \\
                        \frac{\dd}{\dd z}\vecc(z)
                \end{pmatrix}$,
        and rewrite the equation as a first-order linear ODE:
        \begin{equation}
                \frac{\dd}{\dd z} \vecy(z) = B(z)\vecy(z),\quad
                B(z) \coloneqq \begin{pmatrix}
                        O_n & I_n \\
                        -\left(U(z)+\Gamma^2\right) & O_n
                \end{pmatrix}.
                \label{eq:ode2}
        \end{equation}
        Note that following Eqs.~\eqref{eq:vtau} and~\eqref{eq:vrho}, there is a linear transformation $S$ such that:
        \begin{equation}
                \begin{pmatrix}
                        \vtau(z) \\
                        \vrho(z)
                \end{pmatrix}
                =
                S
                \vecy(z), \quad
                S\coloneqq\begin{pmatrix}
                                \Gamma & -\imag I_n \\
                                \Gamma & \imag I_n
                \end{pmatrix}.
        \end{equation}
        The solution of Eq.~\eqref{eq:ode2} has linear form:
        \begin{equation}
                \vecy(z) = Y(z)y(\zbot),
        \end{equation}
        where $Y$ is a $2n\times 2n$ matrix function that solves an initial value problem of a matrix ODE:
        \begin{equation}
                \frac{\dd}{\dd z}Y(z) = B(z)Y(z), \quad Y(\zbot) = I_{2n}.
        \end{equation}
        Once the numerical solution of $Y(0)$, $\hat{Y}_L$ is computed, the desired value $\hat{\vrho}(0)$ is obtained by simple linear algebra via a block approach:
        \begin{equation}
                \begin{pmatrix}
                        \hat{\vtau}(0) \\
                        \hat{\vrho}(0)
                \end{pmatrix}
                = S\hat{Y}_LS^{-1}
                \begin{pmatrix}
                        \vtau(\zbot) \\
                        \vzero
                \end{pmatrix}
                = S\hat{Y}_LS^{-1}\begin{pmatrix}
                        I_n \\
                        O_n
                \end{pmatrix}
                \vtau(\zbot).
                \label{eq:mateq}
        \end{equation}

        The last equality shows that the initial value always lies in the $n$-dimensional subspace of
        $\mathbb{C}^{2n}$ spanned by the first $n$ columns of $S^{-1}$.
        Therefore, we obtain an \emph{economical form} of the ODE Eq.~\eqref{eq:ode2} by changing the initial value of the equation:
        \begin{equation}
                \frac{\dd}{\dd z}Z(z) = B(z)Z(z),\quad Z(\zbot) = S^{-1}\begin{pmatrix}I_n \\ O_n \end{pmatrix}.
        \end{equation}
        The solution of this form is a $2n\times n$ matrix function $Z(z) = Y(z)S^{-1}\begin{pmatrix}I_n \\ O_n\end{pmatrix}$.
        If we let $\hat{Z}_L$ be the numerical solution of $Z(0)$,
        we can rewrite Eq.~\eqref{eq:mateq} for this form as follows:
        \begin{equation}
                \begin{pmatrix}
                        \hat{\vtau}(0) \\
                        \hat{\vrho}(0)
                \end{pmatrix}
                = S\hat{Z}_L\vtau(\zbot).
                \label{eq:mateq2}
        \end{equation}
        This can also be solved by simple linear algebra via a block approach,
        and the numerical solution $\hat{\vrho}(0)$ can be written explicitly as follows.
        Let the upper and lower halves of $\hat{Z}_L$ be $n\times n$ matrices $\hat{Q}_L$ and $\hat{P}_L$, respectively,
        i.e., $\hat{Z}_L = \begin{pmatrix} \hat{Q}_L \\ \hat{P}_L\end{pmatrix}$; then, 
        we have
        \begin{equation}
                \hat{\vrho}(0) = \left(\Gamma+\imag \hat{P}_L\hat{Q}_L^{-1}\right)\left(\Gamma-\imag \hat{P}_L\hat{Q}_L^{-1}\right)^{-1} \vtau(0). \label{eq:vrhocomp}
        \end{equation}
        Because the economical form has half the number of columns as the full-matrix form, the computational cost is roughly halved as well.
        Thus, the economical form is preferable for computation.

        \subsection{Right-hand side transformation}
        The new forms of the ODE also suffer from numerical breakdown
        because they are simply linear transformations of the original ODE.
        As integration proceeds from $z=\zbot$ to $0$, some components of $Z(z)$ grow and others decay; thus,
        $\hat{Q}_L$ becomes numerically singular.
        To avoid this, we developed a new technique called the \emph{right-hand side transformation} (RHST)
        which applies a linear transformation
        to the intermediate state in the ODE solver from the right-hand side.

        Let us consider using a Runge-Kutta-like method to compute the numerical solution for the ODE Eq.~\eqref{eq:mateq2}.
        Runge-Kutta-like methods generate intermediate states $\hat{Z}_i$ for every $z_i$ steps, which approximate the solution of the ODE
        $Z(z_i) \approx \hat{Z}_i$.
        The computation of the RHST is simple: we apply a linear transformation $T_i$ to $\hat{Z}_i$ from the right-hand side
        to improve the numerical condition.
        For $\hat{Z}_0^{(0)} = \hat{Z}_0$, the RHST consists of the five steps listed below:
        \begin{enumerate}
                \item \label{rp1} Compute $\hat{Z}_{i+1}^{(i)}$ from $\hat{Z}_{i}^{(i)}$ using an ODE solver.
                \item \label{rp2} Compute $\xi_{i+1}$, the \emph{badness} of $\hat{Z}_{i+1}^{(i)}$.
                \item \label{rp3} If $\xi_{i+1}$ is greater than the threshold, calculate $T_{i+1}$, if not, let $T_{i+1}=I_n$.
                \item \label{rp4} Compute $\hat{Z}_{i+1}^{(i+1)} = \hat{Z}_{i+1}^{(i)} T_{i+1}$.
                \item Repeat the above process from $i=0$ to $L-1$.
        \end{enumerate}

        The key idea is that linear transformations from the right do not change the final result $\hat{\vrho}(0)$
        because they do not change the subspace spanned by the column of the initial matrix $\hat{Z}_0$.
        Because Runge-Kutta-like methods are linear for linear ODEs, there exist $2n\times 2n$ linear transformations $\bar{Z}_i$, which satisfy the following condition:
        \begin{equation}
                \hat{Z}_i = \bar{Z}_i \hat{Z}_0.
        \end{equation}
        Therefore, the application of $T_{i+1}$ at step \ref{rp3} simply changes the initial value of the ODE:
        \begin{multline}
                \hat{Z}_{i+1}^{(i+1)} = \hat{Z}_{i+1}^{(i)}T_{i+1} = \hat{Z}_{i+1}T^{(i+1)} = \bar{Z}_{i+1} \left(\hat{Z}_0 T^{(i)}\right),\\
                T^{(i)}\coloneqq T_1T_2\cdots T_{i}.
        \end{multline}
        Let us suppose that $T^{(L)}$ is non-singular. Substituting $\hat{Q}_L$ and $\hat{P}_L$
        with $\hat{Q}_LT^{(L)}$ and $\hat{P}_LT^{(L)}$, respectively,
        in Eq.~\eqref{eq:vrhocomp} gives
        \begin{align}
                \hat{\vrho}(0) &=
                \left(\Gamma+\imag \hat{P}_LT^{(L)}(\hat{Q}_LT^{(L)})^{-1}\right)\times {} \nonumber \\
                &{} \quad\quad \left(\Gamma-\imag \hat{P}_LT^{(L)}(\hat{Q}_LT^{(L)})^{-1}\right)^{-1} \vtau(0) \\
                &=
                \left(\Gamma+\imag \hat{P}_L\hat{Q}_L^{-1}\right)\left(\Gamma-\imag \hat{P}_L\hat{Q}_L^{-1}\right)^{-1} \vtau(0).
        \end{align}
        Therefore, the RHST does not change the final result.

        \begin{algorithm}[t]
                \caption{Algorithm of the proposed method with the RHST.} \label{alg:alg1}
                \begin{algorithmic}
                        \Procedure{SolveODE\_STR}{$n$, $L$, $h$, $\Gamma$, $\mathrm{ODE}$}\\
                        \begin{tcolorbox}[size=fbox,colback=white, width=\linewidth,
                                left skip=5mm,before skip balanced=0mm]
                                Compute reflection wave at the surface $\hat{\vrho}(0)$.
                                $n$ and $L$ are the numbers of reciprocal rods and slices, respectively.
                                $h$ is  the step size.
                                $\Gamma$ is a diagonal matrix,
                                and $\mathrm{ODE}(z_{i-1}, z_i, \hat{Z}_{i-1})$ is an ODE solver
                                for a slice
                                that computes the next state $\hat{Z}_i$
                                from the previous state $\hat{Z}_{i-1}$.
                        \end{tcolorbox}
                        \State $\hat{Z}\gets [\Gamma, -\imag I_n; \Gamma, \imag I_n] \backslash [I_n; O_n]$
                        \For{$i=1$ to $L$}
                        \State $\hat{Z} \gets \mathrm{ODE}((i-1)h, ih, \hat{Z})$
                        \State $\xi \gets \operatorname{estcondGC}(\hat{Z}(1\sd n, \sd))$ \Comment{Estimate cond.\ \#}
                        \If{$\xi > 1,000$} \Comment{$1,000$ is a threshold constant}
                        \State $\hat{Z} \gets \hat{Z} / \hat{Z}(1\sd n, \sd)$ \Comment{RHST}
                        \EndIf
                        \EndFor
                        \State $\hat{Z} \gets [\Gamma, -\imag I_n; \Gamma, \imag I_n] * \hat{Z}$
                        \State $\hat{R}_L \gets \hat{Z}((n+1)\sd (2n),\sd) / \hat{Z}(1\sd n, \sd)$
                        \State $\hat{\vrho}(0) \gets \hat{R}_L(\sd,1)$ \Comment{Reflection of the plane wave $\veck_1=\vzero$.}
                        \EndProcedure
                \end{algorithmic}
        \end{algorithm}

        Thus far, we have not discussed how to determine $T_{i+1}$ and $\xi_{i+1}$.
        Let $\hat{Q}_{i+1}^{(i)}$ and $\hat{P}_{i+1}^{(i)}$ be the upper- and lower-half part of
        $\hat{Z}_{i+1}^{(i)}$, respectively.
        We then use the following definitions:
        \begin{itemize}
                \item $T_{i+1}$: inverse of $\hat{Q}_{i+1}^{(i)}$
                \item $\xi_{i+1}$: estimated condition number of $\hat{Q}_{i+1}^{(i)}$ based on the Gershgorin circle theorem
        \end{itemize}
        These definitions keep $\hat{Q}_{i+1}^{(i)}$ close to the identity;
        thus, there is a lower likelihood of it becoming numerically singular.
        $\xi_{i+1}$ has three advantages:
        the computational cost is small,
        overestimation is guaranteed,
        and the accuracy is reasonable because
        $\hat{Q}_{i+1}^{(i)}$ is close to the identity.
        Moreover, these definitions are related to the conventional method,
        as described in the next subsection.

        The algorithm~\ref{alg:alg1} describes the overall process.
        We use Matlab-like syntax for matrix compositions and matrix operations.
        As can be seen in the algorithm, the RHST is simple, adding only four lines of code
        that compute the estimated condition number, compare it with a threshold constant,
        apply the inverse from the right, and \texttt{end if}.
        Therefore, the algorithm is generic to ODE solvers as long as it is a single-step method.

        \subsection{Relationship with the conventional method}
        Our proposed method can be seen as a generalization of the conventional method in three aspects.

        First, the new ODE Eq.~\eqref{eq:ode2} is a linear transformation of the conventional ODE Eq.~\eqref{eq:tr1},
        thus, the equations are essentially the same as a BVP in theory.

        Second, because we explicitly write the ODE of the operator $Z(z)$ in Eq.~\eqref{eq:mateq2},
        our proposed method can be used with a variety of ODE solvers.
        In contrast, the conventional method is tied to a single integration scheme, the second-order Magnus method.

        Third, the RHST can be seen as a generalization of the recursive reflection technique.
        In fact, we can assume that $\hat{Z}_i^{(i)}$ has the form of
        $\hat{Z}_{i}^{(i)} = \begin{pmatrix} I_n \\ \hat{R}_i\end{pmatrix}$
        by choosing $T_i$ appropriately.
        Now let us consider using the second-order Magnus method for the proposed method.
        Here we apply $\hat{A}_i=\ee^{(z_{i+1}-z_i)A_i}$ to $\hat{Z}_{i}^{(i)}$ from the left:
        \begin{equation}
                \hat{Z}_{i+1}^{(i)} = \hat{A}_i \hat{Z}_{i}^{(i)} = \begin{pmatrix}\hat{X}_i+\hat{Y}_i\hat{R}_i \\ \hat{Z}_i+\hat{W}_i\hat{R}_i\end{pmatrix}.
        \end{equation}
        Then, by applying $T_i=\left(\hat{X}_i+\hat{Y}_i\hat{R}_i\right)^{-1}$ from the right, we find
        \begin{equation}
                \hat{Z}_{i+1}^{(i+1)}=\hat{Z}_{i+1}^{(i)}T_i=\begin{pmatrix}
                        I_n \\
                        \left(\hat{Z}_i+\hat{W}_i\hat{R}_i\right)\left(\hat{X}_i+\hat{Y}_i\hat{R}_i\right)^{-1}
                \end{pmatrix}.
        \end{equation}
        The lower-half part of this matrix has the same form as Eq.~\eqref{eq:ri1},
        thus, by letting $\hat{R}_{i+1}$ be that part of the matrix,
        we reproduce $\hat{R}_i$ for $i=0,\ldots,L-1$ in the conventional method using the RHST.
        Note that we must change the initial value from $\hat{Z}_0$ to $S\hat{Z}_0$.

        \subsection{Choice of ODE solver}
        Here, we describe how two types of concrete ODE solvers,
        explicit Runge-Kutta methods and the splitting methods, can be used for the proposed method
        and compare their computational patterns with that of the conventional method.

        The $s$-step Runge-Kutta methods can be described by the coefficients $a_{r,j}$, weights $b_{j}$, and nodes $c_r$ where $1\le r, j\le s$.
        Let $h=z_{i+1}-z_i$, $p_0=\hat{P}_{i}^{(i)}$ and $q_0=\hat{Q}_{i}^{(i)}$ to simplify the notation.
        When applying the Runge-Kutta method to the proposed method, the computation consists of the following steps:
        \begin{align}
                p_r &= q_0 + h\sum_{j=1}^{r-1}a_{r,j}q_j, \\
                q_r &= -\left(U(z_i+hc_r)+\Gamma^2\right)\left(p_0 + h\sum_{j=1}^{r-1}a_{r,j}p_j\right),
        \end{align}
        and the final result is calculated as follows:
        \begin{align}
                \hat{P}_{i+1}^{(i)} &= p_0 + h\sum_{j=1}^sb_jp_j, \\
                \hat{Q}_{i+1}^{(i)} &= q_0 + h\sum_{j=1}^sb_jq_j.
        \end{align}

        The splitting method is known as a type of geometric integrator~\cite{Blanes2002},
        but also serves as a simple and storage-efficient variant of the Runge-Kutta-Nystr\"{o}m method.
        Now, let $\tau_r$ be a pseudo time variable, where $\tau_0=z_i$, $a_r$, and $b_r$ are the nodes.
        There are two major variants of the splitting method, known as ``ABA'' and ``BAB.''
        Applying the $s$-step ``BAB'' method to the proposed method gives rise to a computational formula consisting of the following three steps for $r=1$ to $s$:
        \begin{align}
                q_r &= q_{r-1} - hb_r\left(U(\tau_r)+\Gamma^2\right)p_{r-1},\\
                p_r &= p_{r-1} + ha_r q_r, \\
                \tau_r &= \tau_{r-1} + ha_r,
        \end{align}
        and the final result is computed as follows:
        \begin{align}
                \hat{Q}_{i+1}^{(i)} &= q_s - hb_s\left(U(\tau_s)+\Gamma^2\right)p_s,\\
                \hat{P}_{i+1}^{(i)} &= p_s.
        \end{align}
        The ``ABA'' method is similar but alternates the roles of $p_r$ and $q_r$.

        We note that both methods are linear because each step can be written as linear transformations from the left.
        Additionally, both methods can be used with the RHST because they are single-step methods.

        Both methods are similar in computation pattern: they consist of
        matrix multiplications of $n\times n$ square matrices and a weighted sum of matrices for each intermediate step.
        This approach is far more efficient and requires less computation than
        the matrix exponential computation in the conventional method.
        The Runge-Kutta method requires storage to hold all of the intermediate states $p_r$ and $q_r$,
        while the splitting method does not which is the same as the conventional one.

        One of the main disadvantages of these ODE solvers arises because
        the number of evaluations of the potential $U(z)$ is multiplied by $s$.
        This may cause a performance problem if the computation time of $U(z)$ is large.
        In RHEED\slash TRHEPD simulations, the computational cost of $U(z)$ is not large
        because the same $U(z)$ can be used for simulations of different angles.
        Moreover, this problem is mitigated by the reduction in the number of slices $L$ by using higher-order ODE solvers.

        \subsection{Improvements in implementation}\label{ss:hpc}
        We also improved the implementation of the original \texttt{sim-trhepd-rheed} code for the optimal performance on recent CPUs. In particular, 
        linear algebraic procedures are reimplemented by the packages like
        the basic linear algebra subroutines (BLAS)~\cite{Lawson1979} and the linear algebra package (LAPACK)~\cite{lapack1999}, since
        most of the computation for the conventional method and proposed method consists of matrix computations.
        The original code uses LAPACK's subroutine for eigenvalue decomposition,
        but homemade subroutines for other matrix computations.
        We replaced these subroutines with functionally equivalent LAPACK's subrouitnes,
        and reordered \texttt{do}-loops to split out matrix-matrix multiplications and 
        replaced them with the BLAS's general matrix-matrix multiplication subroutine, \texttt{GEMM}.
        The use of these libraries enable us the optimal computation among the recent CPUs such as SIMD and the multi-level cache-memory hierarchy.

        Another improvement of this implementation arises from multi-threading based on the parallelism of
        the angles of the incident wave.
        RHEED\slash TRHEPD measurements are performed for many pairs of angles $(\theta_0, \theta_1)$,
        and the simulation for each pair is trivially parallelized.
        The number of pairs is sufficiently large for multi-threading, $\sim 100$;
        thus, we added OpenMP~\cite{OpenMP2018} directives before the \texttt{do}-loop in the source code for thread parallelization.
        The original code has no explicit parallelization, and most of the code runs on a single core of a CPU.

        We made a few minor changes such as using real arithmetic as much as possible
        in the potential computation to reduce the number of computations and
        increasing the number of digits for output from $4$ to $15$ for error analysis, as discussed in the next section.

        To observe the effect of these improvements, we developed a new implementation of the conventional method named \texttt{opt}
        which is compared with the original implementation \texttt{orig}.
        The implementations of the proposed method discussed in the next section include the improvements explained in this subsection.

        \section{Performance evaluation}
        In this section, we evaluate the performance of the proposed method using time-error charts.
        We use two definitions of the error:
        \texttt{eorig} is the difference from the result of the reference (default) settings of the conventional method, and
        \texttt{eacc} is the difference from the ``accurate'' result of the proposed method with a fine step size.
        Let $\hat{\veta}$, $\hat{\veta}_\mathrm{acc}$ and $\hat{\veta}_\mathrm{orig}$
        be the computation results obtained by the target implementation, the new implementation with the fine step size,
        and the conventional implementation with the default settings, respectively.
        Then, we define the two errors as follows:
        \begin{align}
                \mathtt{eorig} &:= \frac{\max_{\vecb_0 \in B_0}\left\|\hat{\veta}(\vecb_0)-\veta_\mathrm{orig}(\vecb_0)\right\|}{\max_{\vecb_0 \in B_0} \left\|\hat{\veta}_\mathrm{orig}(\vecb_0)\right\|}, \\
                \mathtt{eacc} &:= \frac{\max_{\vecb_0 \in B_0}\left\|\hat{\veta}(\vecb_0)-\veta_\mathrm{acc}(\vecb_0)\right\|}{\max_{\vecb_0 \in B_0} \left\|\hat{\veta}_\mathrm{acc}(\vecb_0)\right\|}.
        \end{align}
        Here, $B_0$ is the set of projected wavevectors of the incident wave.
        We used \texttt{eacc} instead of the error based on the theoretical solution
        because the latter is difficult to compute for this simulation.
        We also provide \texttt{eorig} to avoid the case
        in which the results of the same method cause unexpected relationships.
        \texttt{eorig} is also useful for researchers who are using the conventional implementation.

        We consider the following five combinations of methods and implementations:
        \begin{itemize}
                \item \texttt{orig}: The conventional implementation of the conventional method~\cite{Hanada2022}, cloned the commit \texttt{df61124c} from \cite{str-github} and built without modifications except for the \texttt{Makefile}.
                \item \texttt{opt}: A modified version of \texttt{orig} using the improvements in \S~\ref{ss:hpc}.
                \item \texttt{rk4}: An implementation of the proposed method with the fourth-order Runge-Kutta method, using the improvements in \S~\ref{ss:hpc}.
                \item \texttt{sp4} and \texttt{sp6}: implementations of the proposed method with the fourth- and sixth-order splitting method ($\mathrm{SRKN}_6^b$ and $\mathrm{SRKN}_{11}^b$, respectively, in ~\cite{Blanes2002}),
                using the improvements in \S~\ref{ss:hpc}.
        \end{itemize}
        In the new implementations, we use a step size \texttt{dz} $\approx h$ close to the E-6 series preferred numbers from $0.01\si{\angstrom}$ to $0.69\si{\angstrom}$.
        \texttt{orig} has a minor bug that causes instability of the computational domain when \texttt{dz} changes.
        Thus, we fix $\mathtt{dz}=0.01\si{\angstrom}$ for \texttt{orig}.

        \begin{table}
                \centering
                \caption{Details of the test data used in the experiments}\label{tab:testdata}
                \begin{tabular}{lccc}
                        \toprule
                        & \texttt{n=23} & \texttt{n=47} & \texttt{n=521} \\
                        \midrule
                        surface structure & \multicolumn{3}{c}{Si 7x7 (111)} \\
                        cell symmetry & \multicolumn{3}{c}{p3m1} \\
                        \# of atoms in a cell & \multicolumn{3}{c}{37} \\
                        domain size w/o bulk layer& \multicolumn{3}{c}{9.910955\si{\angstrom}}\\
                        \# of reciprocal rods $n$ & $23$ & $47$ & $521$ \\
                        \# of angles & \multicolumn{2}{c}{$69$} & $1$ \\
                        altitude $\theta_0$ & \multicolumn{2}{c}{$0.1^{\circ}, 0.2^{\circ}, \ldots 69.0^{\circ}$}& $1.3^{\circ}$ \\
                        azimuth $\theta_1$ & \multicolumn{2}{c}{$60^{\circ}$} & $-30^{\circ}$ \\
                        \bottomrule
                \end{tabular}
        \end{table}
        We used three test datasets, denoted as \texttt{n=23}, \texttt{n=47}, and \texttt{n=521}.
        Details are listed in Table~\ref{tab:testdata}.
        The numerical problem is one for the TRHEPD simulator for the Si(111)-$7 \times 7$ surface, a famous semiconductor surface. The atom positions of the Si(111)-$7 \times 7$ surface and their RHEED and TRHEPD diffraction images are found in Ref.~\cite{Fukaya2019a}. 
        The main difference among the test datasets is the number of reciprocal rods $n$.
        The number of glancing angles for \texttt{n=23} and \texttt{n=47} is  $69$ while that of \texttt{n=521} is $1$,
        because the calculation for a higher number of angles is too time-consuming for \texttt{orig}.
        Therefore, multi-threading based on the number of angles cannot be applied to \texttt{n=521},
        instead, the BLAS and LAPACK implementations parallelize the matrix computations.

        All of the time measurements were performed on a BTO desktop PC with an Intel i7-12700 processor (uses 8 P-cores, fixed to 2.1 GHz) and dual channel DDR4-3200 memories.
        We used the Intel oneMKL version~\texttt{2022.1.0} ~\cite{intelmkl2023} for the implementation of BLAS and LAPACK,
        which is highly optimized for Intel CPUs.

        \paragraph{Time--error chart}
        \begin{figure*}[t]
                \centering
                \includegraphics[width=0.7\linewidth]{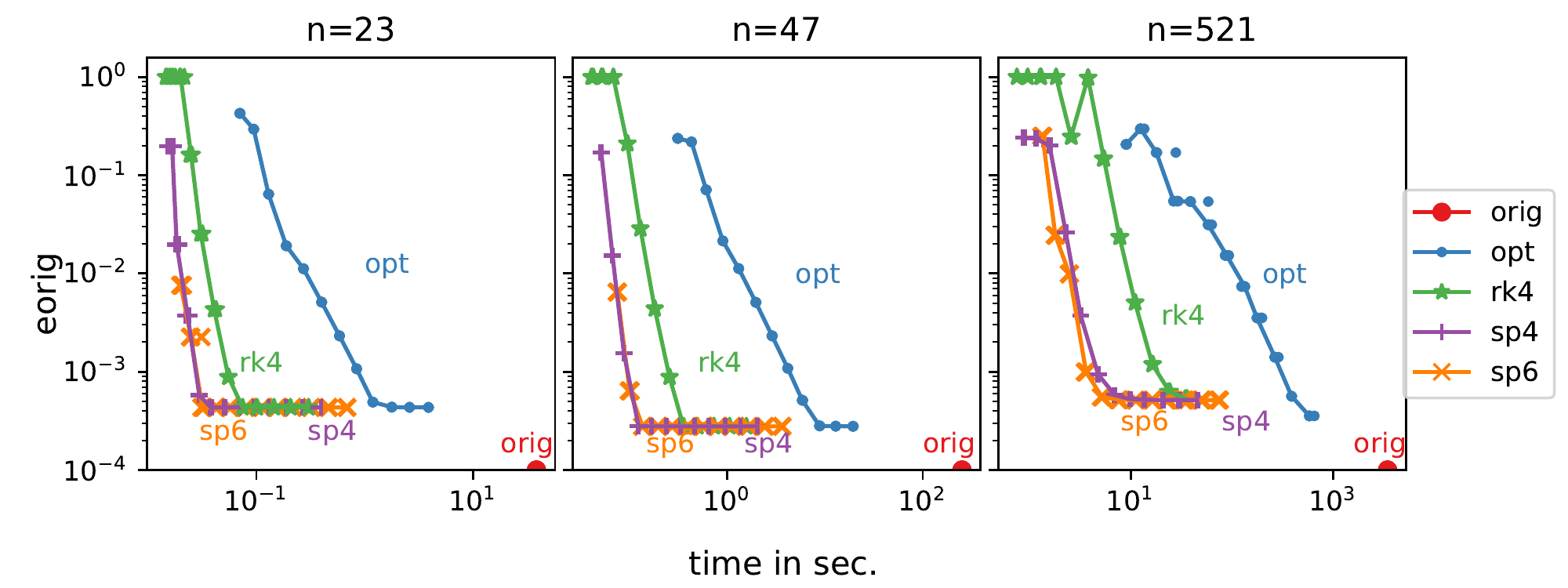}
                \includegraphics[width=0.7\linewidth]{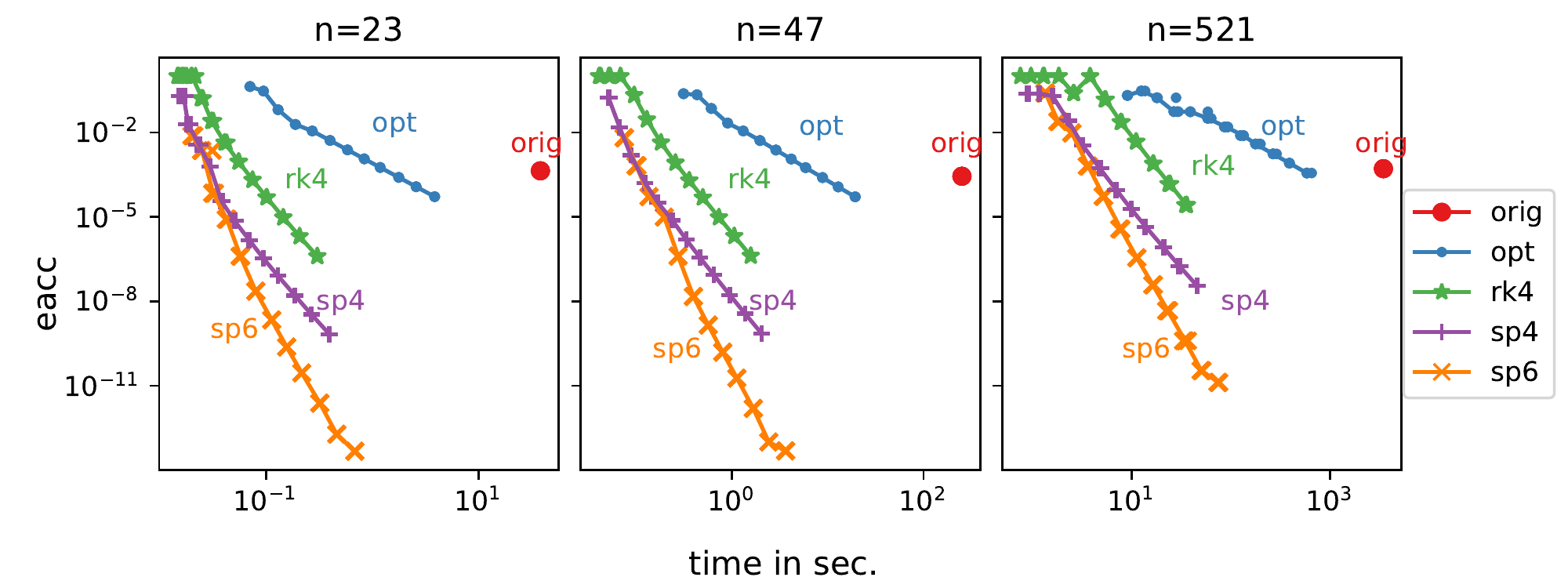}
                \caption{Time--error charts for the conventional and proposed implementations.
                        The top figures are for \texttt{eorig}, and the bottom figures are for \texttt{eacc}.
                        All measurements were performed five times, and log-log-plotted as markers.
                        The lines are the medians for each \texttt{dz}.
                        We placed the results for \texttt{orig} on the x-axis in the top figures
                        to indicate that the \texttt{eorig} value is exactly $0$ by definition.}
                \label{fig:time-error}
        \end{figure*}
        Time--error charts are shown in Fig.~\ref{fig:time-error}.
        The curves of \texttt{eorig} are saturated around $\mathtt{eorig} \sim 10^{-4}$.
        This occurs simply because the output format of \texttt{orig} is \texttt{'E12.4'}, only $4$ digits.
        Consequently, the new implementations achieve an accuracy that is considered to be sufficient by the developers of the conventional implementation more quickly than \texttt{orig}.
        Among the new implementations, \texttt{sp4} and \texttt{sp6} achieve a sufficient accuracy within the shortest amount of time.
        The performances of \texttt{sp4} and \texttt{sp6} appear to be almost the same in this figure,  even though \texttt{sp6} has a higher order than \texttt{sp4}.

        In the figures for \texttt{eacc}, we can see the behaviors of the curves at values below $\mathtt{eorig} < 10^{-4}$.
        The curves for \texttt{sp6} are the steepest among the new implementations
        and outperform \texttt{sp4} at values below $\mathtt{eacc}<10^{-6}$ for \texttt{n=23} and \texttt{n=47}
        and at values below $\mathtt{eacc}<10^{-4}$ for \texttt{n=521}.
        \texttt{sp4} and \texttt{rk4} have almost the same slope, but \texttt{sp4} is faster than \texttt{rk4} by one order of magnitude.

        \paragraph{Computation time for baseline accuracy}
        \begin{table*}
                \centering
                \caption{Computation time and speed-up rate from \texttt{orig}.
                        We selected the fastest results such that the $\mathtt{eacc}$ value
                        is smaller than that of \texttt{orig}.}\label{tab:times}
\begin{tabular}%
        {%
                l%
                l%
                S[exponent-mode = fixed, round-mode = places, round-precision=3]%
                S[exponent-mode = scientific, print-zero-exponent = true, round-mode = figures, round-precision=3, table-alignment=center]%
                S[exponent-mode = fixed, round-mode = figures, round-precision=3]%
                S[exponent-mode = fixed, round-mode = figures, round-precision=3]%
        }
        \toprule
        data & method & {\texttt{dz}} & {$\mathtt{eacc}$} &  {time in sec.} &  {speed up rate} \\
        \midrule
        n=23 &   \texttt{orig} & 1.000000000000000e-02 & \num{4.329095380798452e-04} & 3.830600000000000e+01 & 1.000000000000000e+00 \\
        &    \texttt{opt} & 2.200000000000000e-02 & \num{2.524585847350836e-04} & 1.775955341800000e+00 & 2.156923606039293e+01 \\
        &    \texttt{rk4} & 4.700000000000000e-02 & \num{2.046810891250196e-04} & 7.463165259999868e-02 & 5.132674765398493e+02 \\
        &    \texttt{sp4} & 1.500000000000000e-01 & \num{3.562029263491765e-05} & 3.813462160000256e-02 & 1.004494036988096e+03 \\
        &    \texttt{sp6} & 3.300000000000000e-01 & \num{6.913115388572034e-05} & 3.186411620000058e-02 & 1.202167345849665e+03 \\
        \midrule
        n=47 &   \texttt{orig} & 1.000000000000000e-02 & \num{2.777707069907318e-04} & 2.516040000000000e+02 & 1.000000000000000e+00 \\
        &    \texttt{opt} & 2.200000000000000e-02 & \num{2.490417808814971e-04} & 8.813310090600044e+00 & 2.854818421382355e+01 \\
        &    \texttt{rk4} & 4.700000000000000e-02 & \num{1.962351123813380e-04} & 3.579840540000305e-01 & 7.028357749140931e+02 \\
        &    \texttt{sp4} & 2.200000000000000e-01 & \num{1.597960656229034e-04} & 1.232810935999851e-01 & 2.040896885749482e+03 \\
        &    \texttt{sp6} & 3.300000000000000e-01 & \num{5.223740308576600e-05} & 1.361412525999867e-01 & 1.848109924030657e+03 \\
        \midrule
        n=521 &  \texttt{orig} & 1.000000000000000e-02 & \num{5.151892346217075e-04} & 3.441453999999999e+03 & 1.000000000000000e+00 \\
        &   \texttt{opt} & 1.000000000000000e-02 & \num{3.611137709371056e-04} & 5.955675777463991e+02 & 5.778444174248548e+00 \\
        &   \texttt{rk4} & 1.500000000000000e-02 & \num{1.470262355964512e-04} & 2.396980442839995e+01 & 1.435745548229208e+02 \\
        &   \texttt{sp4} & 6.800000000000000e-02 & \num{9.103368786673352e-05} & 6.966074966000815e+00 & 4.940305719930715e+02 \\
        &   \texttt{sp6} & 1.500000000000000e-01 & \num{5.381183283290355e-05} & 5.183526088799772e+00 & 6.639214197139030e+02 \\
        \bottomrule

\end{tabular}

        \end{table*}

        Table~\ref{tab:times} shows the fastest results from Fig.~\ref{fig:time-error},
        where the $\mathtt{eacc}$ value is less than that of \texttt{orig} with $\mathtt{dz}=0.01\si{\angstrom}$.
        This table clearly shows that the new implementations outperform \texttt{orig}.
        \texttt{opt} is more than $20$ times faster than \texttt{orig} then the two improvements described in \S \ref{ss:hpc} are implemented
        for \texttt{n=23} and \texttt{n=47} and approximately $6$ times faster when one of the two improvements is included for \texttt{n=521}.
        With the proposed method, \texttt{rk4}, \texttt{sp4}, and \texttt{sp6} are more than $100$ times faster than \texttt{orig},
        and the performance of \texttt{sp4} shows as improvement of more than $2,000$ fold for \texttt{n=47}.

        \section{Discussion}
        \paragraph{Effect of the recursive reflection technique and RHST}
        It is not fully theoretically understood why the recursive reflection technique and RHST can avoid numerical breakdowns.
        One reason might be that the physical law bounds some norm of the transfer matrix $R_i$ less than $1$ because the reflection wave must have less energy than the input.

        The transfer matrix $\hat{A}_i$ should have the same property as a physical representation.
        However, the transfer matrix is an operator which converts the wave at the bottom side to that of the top side
        and not an operator which converts the input to the output.
        This is the reason why $\hat{A}_i$ can have a norm greater than $1$.
        If we convert the transfer matrix to separate the input (from the top and bottom) and output (to the bottom and top)
        at the left and right side of the matrix,
        the norm of the converted matrix will be less than or equal to $1$.

        \paragraph{Need for structure preservation}
        The conventional method uses the second-order Magnus method ~\cite{Blanes2008}, which is a Lie-group method that can maintain a Lie-group structure if the equation has such a structure.
        The ODE solvers that we used for our proposed method in the test cannot maintain such a structure.
        Thus, this might be a regression from the conventional method.
        
        As far as we know, the equation has no interesting Lie-group structure because the potential $v$ has an artificial imaginary component that represents the absorption effect.
        Therefore, the structure is $\mathrm{GL}(2n)$, which is the group of non-singular matrices.

        If the geometry is important,
        higher-order Magnus-based methods for second-order non-autonomous ODEs are available~\cite{bader2018}.

        \paragraph{Application of the proposed method to other simulations}
        Our goal in this study was to improve upon conventional methods and implementations currently in use.
        At this stage, the proposed method is specific to RHEED\slash TRHEPD,
        but may possibly be applicable to other many-beam reflective diffraction simulations
        that have severe ill-conditionedness.

        Note that if the problem is well-conditioned, i.e.,
        if the matrix product $\hat{A}_0\hat{A}_1\cdots$ has a small condition number,
        we can use simple linear solvers to solve BVP.
        The generalized minimum residual method (GMRES) will be the best method in this case,
        because the method does not require computing the explicit components of the matrix product.

        \section{Conclusion}
        In this article, we proposed a new method for faster RHEED\slash TRHEPD simulations 
        of the conventional ones.
        Our strategy is standard, reformulates BVP as an initial-value matrix ODE,
        and applies high-order ODE solvers such as fourth-order Runge-Kutta and splitting methods,
        except for the generalization of the \emph{recursive reflection technique} to the \emph{RHST}.
        As a result, our proposed method reduces the number of computations and
        increases the computation efficiency while maintaining the same accuracy as that of the conventional method.
        Moreover, we also proposed a high-performance implementation of the algorithm to
        utilize the multi-thread parallelism and cache memory performance of recent CPUs.

        In our performance evaluation based on three test problems,
        our new implementations of the proposed algorithm outperform the conventional method
        by orders of magnitude, up to $2,000$ fold.
        This huge leap in speed from the conventional method not only reduces the cost of simulations,
        but also widens the ability of reverse analysis, for example, by allowing the number of parameters to be increased for optimization
        or enabling reverse analysis to be performed immediately after measurements without supercomputers.

        Most of the computation of the algorithm presented herein consists of  matrix multiplications, which are well-suited for accelerators such as graphic processing units (GPUs), except for the RHST step.
        As a future work, it may be beneficial to develop a better method for the RHST that emplys uniform computational patterns with lower or similar computational costs while maintaining the ability to avoid numerical breakdowns.
        The accuracy of the simulation with lower-precision (single or half) could also be evaluated in the future
        to utilize the processing power of these accelerators.

        \section*{Declaration of competing interest}
        The authors have no conflicts of interest to declare that are relevant to the content of this article.

        \section*{Data statement}
        All the source codes and obtained data are available at \url{github.com/shuheikudo/trhepd-opt}.

        \section*{Acknowledgements}
        The authors thank Izumi Mochizuki for providing the input files of the present numerical examples. The present research 
        was supported by the Research Institute for Mathematical Sciences,        an International Joint Usage/Research Center located in Kyoto University,        and was partially supported by Japanese KAKENHI projects (20H00581,21K19773, 22H03598). 

        \bibliographystyle{elsarticle-num}
        \bibliography{library-short}
\end{document}